\documentclass{article}
\usepackage{graphicx}
\usepackage{cite}
\usepackage{amsmath}
\usepackage[figuresright]{rotating}
\title{Data assimilation in 2D incompressible Navier-Stokes equations, using a stabilized explicit $O(\Delta \lowercase{t})^2$ leapfrog finite
difference scheme run backward in time.}{
\author{Alfred S. Carasso\thanks
{Applied and Computational Mathematics Division,
National Institute of Standards and Technology,
Gaithersburg, MD 20899. (alfred.carasso@nist.gov).}}
\begin{document}
\maketitle
\begin{abstract} A stabilized backward marching leapfrog finite difference scheme is applied to hypothetical data at time $T > 0$, to retrieve initial values
	at time $t=0$ that can evolve into useful approximations to the data at time $T$, even when that data does not correspond to an actual solution.
Richardson's leapfrog scheme is notoriously unconditionally unstable in well-posed, forward, linear dissipative evolution equations. Remarkably, that scheme
can be stabilized, marched backward in time, and provide useful reconstructions in an interesting but {\em limited class} of ill-posed, time-reversed, 2D
incompressible Navier-Stokes initial value problems.
Stability is achieved by applying a compensating smoothing operator at each time step to quench the instability. Eventually, this leads
to a distortion away from the true solution. This is the {\em stabilization penalty}. In many interesting cases, that penalty is sufficiently small to allow
for useful results. Effective smoothing operators based on $(-\Delta)^p$, with real $p > 2$, can be efficiently
synthesized using FFT algorithms. Similar stabilizing techniques were successfully applied in several
other ill-posed evolution equations. The analysis of numerical stabilty is restricted to a related linear problem. However, as is found in leapfrog computations of well-posed
meteorological and oceanic wave propagation problems, such linear stability is necessary but not sufficient in the presence of nonlinearities. Here, likewise, additional
Robert-Asselin-Williams (RAW) time-domain filtering must be used to prevent characteristic leapfrog nonlinear instabilty, unrelated to ill-posedness.

Several 2D Navier-Stokes backward reconstruction examples are included, based on the {\em stream function-vorticity} formulation, and focusing on $256 \times 256$ pixel images
of recognizable objects.
Such images, associated with non-smooth underlying intensity data, are used to create severely distorted data at time
$T > 0$. Successful backward recovery is shown to be possible at parameter values significantly exceeding expectations.

\end{abstract}

\pagestyle{myheadings}
\thispagestyle{plain}
\markboth{A.~S. CARASSO}{DATA ASSIMILATION IN 2D NAVIER-STOKES EQUATIONS}
\section{Introduction}
This paper discusses the possible use of direct, non iterative, {\em explicit}, backward marching finite difference schemes, in data assimilation problems for the 2D incompressible
Navier-Stokes equations, at high Reynolds numbers. The actual problem considered here is the ill-posed time-reversed problem of locating appropriate initial values,  given {\em hypothetical} data at some later time. That problem is emblematic of much current geophysical research activity, using computationally intensive iterative methods  based on neural networks
coupled with machine learning, \cite{cintra, howard, blum, qizhi, arcucci, Antil, Chong, Lund,auroux1,ou,auroux2,auroux3,pozo,gosse,gomez, xu, tomislava,camposvelho}.
Previous work on backward marching stabilized explicit schemes in 2D nonlinear data assimilation problems, is discussed in \cite{Burg,Advec,Thermoassim}. There, it is emphasized that the given hypothetical data at time $T > 0$ may
not be smooth, 
{\em may not correspond} to an actual solution at time $T$, and may differ from such
a solution by an {\em unknown large $\delta > 0$} in an  apropriate ${\cal{L}}^p$ norm. Moreover, it {\em may not be possible} to 
locate initial values that can evolve into a useful approximation to the desired data at time $T$.

As is well-known, \cite[p.59]{richtmyer}, for ill-posed initial value problems, all consistent stepwise marching schemes, whether explicit or implicit, are necessarily
unconditionally unstable and lead to explosive error growth. Nevertheless, it is possible to stabilize such schemes by applying an appropriate compensating smoothing operator
at each time step to quench the instability. These schemes then become unconditionally stable, but slightly inconsistent, and lead to distortions away from the true solution.
This is the {\em stabilization penalty}. In many cases, that penalty is sufficiently small to allow for useful results. Indeed, in \cite{carGEM, car1IPSE, car2IPSE,car3IPSE,Thermo, car5IPSE,car6IPSE,car7IPSE},
backward marching stabilized explicit schemes have been successfully applied to several challenging ill-posed {\em backward recovery} problems. 

However, such backward recovery problems differ fundamentally from the data assimilation problem considered here. In \cite{carGEM, car1IPSE, car2IPSE,car3IPSE,Thermo, car5IPSE,car6IPSE,car7IPSE},
the given data at time $T > 0$ are  noisy but relatively smooth, and are known to approximate an {\em actual solution} at time $T$, to within a {\em known small $\delta > 0$} in an  apropriate ${\cal{L}}^p$ norm. While that is not the case in the present
data assimilation problem, several computational examples in Section 7 below,  will show that
success can still be achieved, and on time intervals $[0, T]$
that are {\em several orders of magnitude larger} than might be expected, based on the best-known uncertainty estimates for time-reversed Navier-Stokes
equations.

As was the case in \cite{carGEM, car1IPSE, car2IPSE,car3IPSE,Thermo, car5IPSE,car6IPSE,car7IPSE}, the computational examples in Section 7 will
be based on $256 \times 256$ pixels gray-scale images, with intensity values that are integers ranging from $0$ to $255$.
These images represent easily recognizable objects, are defined by highly non-smooth intensity data that would be quite difficult to synthesize mathematically,
and they present significant challenges to any ill-posed reconstruction procedure. In addition, images allow for easy visual recognition of possible success or failure in 
data assimilation.

\section{2D Navier-Stokes equations in the unit square}
Let $\Omega$ be the unit square in $R^2$ with boundary $\partial \Omega$.
Let $<~,~>$ and $\parallel~\parallel_2$, respectively denote the scalar product and norm on  ${\cal{L}}^2(\Omega)$. We study the 2D problem in
{\em stream function-vorticity} formulation \cite{johnston, ghadimi}. Here, the flow is governed by a stream function $\psi(x,y,t)$ which defines the velocity components
$u(x,y,t) \equiv \psi_y(x,y,t),~v(x,y,t) \equiv -\psi_x(x,y,t),$ together with the vorticity $\omega(x,y,t) \equiv (v_x-u_y)(x,y,t) =-(\Delta \psi)(x,y,t)$.
The present limited data assimilation study will focus on dissipative properties of Navier-Stokes systems in the absence of active boundary conditions or forcing terms.
We restrict attention
to stream
functions $\psi(x,y,t)$ defined on the unit square $\Omega$, which, for the short duration $0 \leq t \leq T,$ of the flow, {\em vanish together with all their derivatives on and near the
boundary}  $\partial \Omega$. However, $\psi(x,y,t)$  will generally have large spatial derivatives away from the boundary.
Given such an initial $\psi(x,y,0)$, and a kinematic viscosity $\nu > 0$, the {\em forward} Navier-Stokes initial value problem takes the form
\begin{equation}
\begin{array}{l }
\omega_t+ u \omega_x + v \omega_y -\nu \Delta \omega =0, ~~~~~~~\Delta \psi =- \omega, \qquad (x,y,t) \in \Omega \times (0,T],  \\
\ \\
u(x,y,0)= \psi_y(x,y,0),~~~~~v(x,y,0)=-\psi_x(x,y,0), \\
\ \\
\omega(x,y,0)=-(\Delta \psi)(x,y,0),\\
\ \\
\psi(x,y,t) = u(x,y,t) = v(x,y,t) = \omega(x,y,t)= 0,  \qquad (x,y,t) \in \partial \Omega \times [0,T].  \\
\label{eq:0}
\end{array}
\end{equation}
For each $t > 0$, consider the Poisson problem on $\Omega:~ \Delta \psi =- \omega,~\psi=0 \in \partial \Omega.$
Let $H$ be the solution operator for that problem, so that for $(x,y,t) \in \Omega \times (0,T]$, we have
$\psi(x,y,t)=-(H \omega)(x,y,t).~$ We may then rewrite the evolution equation in Eq.~(\ref{eq:0}) in the following equivalent form
\begin{equation}
\omega_t -(H \omega)_y \omega_x + (H \omega)_x \omega_y -\nu \Delta \omega =0, \qquad (x,y,t) \in \Omega \times (0,T].
\label{eq:0.01}
\end{equation}
Solving the system in Eq.~(\ref{eq:0}) may be visualized as a marching computation as follows. Given the initial values, the evolution equation for $\omega$ allows advancing
$\omega$ forward one time
step $\Delta t > 0$. The Poisson equation $\Delta \psi =- \omega$ is then solved to obtain $-(H \omega)(x,y, \Delta t)$,  which now 
allows advancing
$\omega$ to the next time step through the evolution equation Eq.~(\ref{eq:0.01}), and so on.

With $A$ the area of the flow domain,
define $U_{max}$ and the Reynolds number $RE$ as follows
\begin{equation}
\begin{array}{l}
	U_{max}=\sup_{\Omega \times [0,T]}\{u^2(x,y,t)+v^2(x,y,t) \}^{1/2},  \\
\ \\
RE=(\sqrt{A}/\nu) U_{max}.
\label{eq:0.1}
\end{array}
\end{equation}

Gray-scale image intensity data have integer values ranging between $0$ and $255$. When such images are used to simulate examples of non-smooth stream functions $\psi(x,y)$,
the primary variables $u, v, \omega,$ in Eq.~(\ref{eq:0}) are derived by taking first and second {\em derivatives} of $\psi$.
To prevent very large derivatives and extremely large Reynolds numbers, such intensity values must be {\bf rescaled}.
Premultiplying by $0.0025$, creates an array $\psi(x,y,0)$ of non negative numbers ranging from $0$ to $0.6375$. This leads to values for $U_{max}$ and $RE$ in Eq.~(\ref{eq:0.1}),
that are $400$ times smaller than they would otherwise be.

The following example is typical.

\subsection{An example: the 2004 Hurricane Ivan stream function}
\begin{figure}
        \centerline{\includegraphics[width=5.0in]{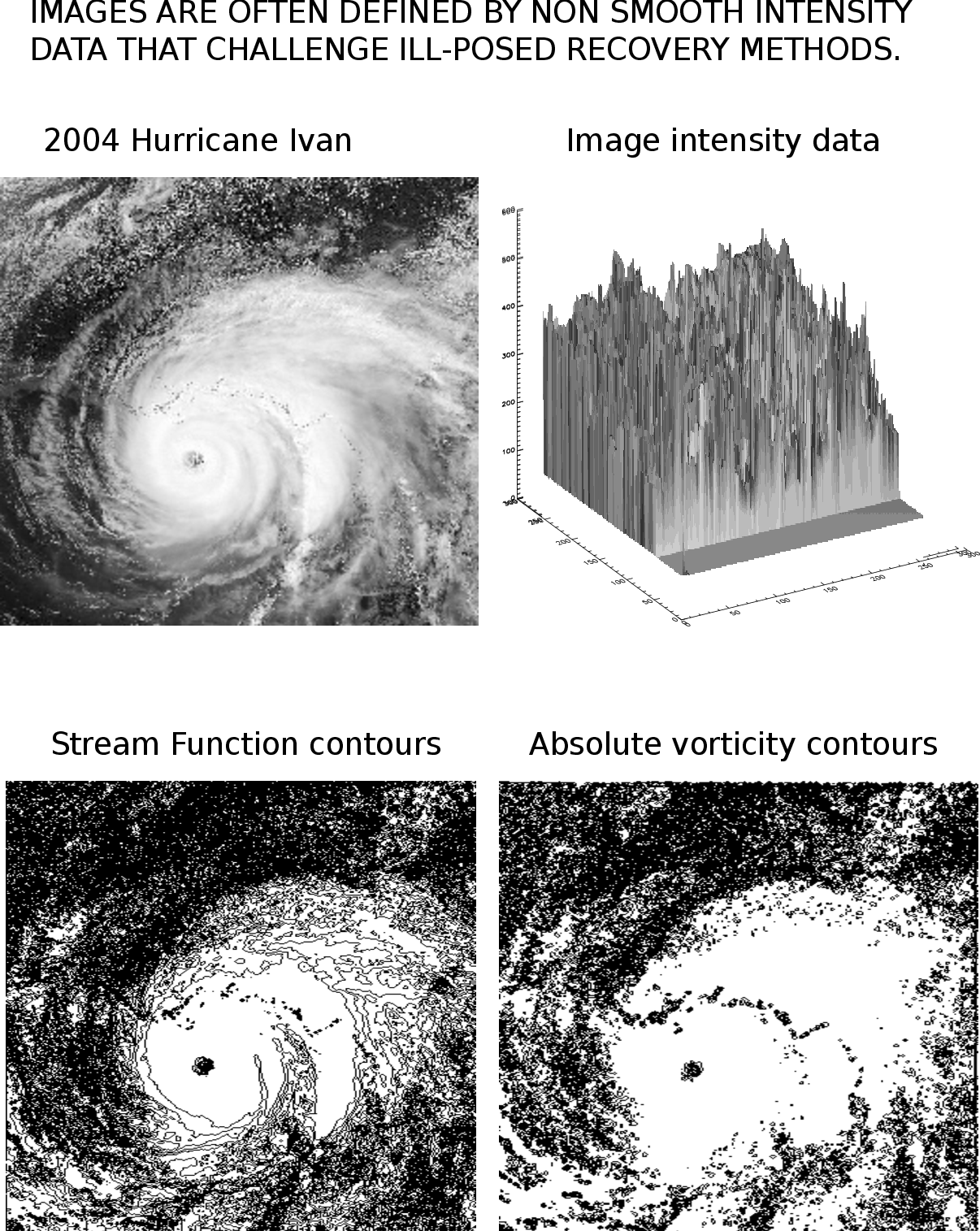}}
\caption{Non-smooth Hurricane Ivan data leads to challenging Navier-Stokes data assimilation problem. Here, with $\nu=0.01,~U_{max}=111, ~RE=11,100,$ and
$\sup_{\Omega}\{|\omega(x,y)|\approx 1.5\times10^5$,
theoretical uncertainty estimates in Eq.~(\ref{eq:3c}) indicate that useful identification of initial values, given the above hypothetical data at some $T > 0$, 
	might not be feasible for  $ T > 1.0 \times 10^{-9}$.}
\end{figure}

The $256 \times 256$ pixel gray-scale Hurricane Ivan image, and associated intensity data, are shown in Figure 1.
The original intensity values were multiplied by $0.0025$, to create an array $\psi(x,y)$ of non negative numbers ranging from $0$ to $0.6375$.
We now view $\psi(x,y)$ as a candidate stream
function, defined on the unit square $\Omega$. Contour plots of $\psi(x,y)$ and of the associated vorticity $\omega(x,y)=-\Delta \psi$,  are also shown in Figure 1.
However, even after reduction by a factor of $400$, one must still expect significant values for the derivatives of $\psi(x,y)$. Indeed, with $U_{max}$ as in Eq.~(\ref{eq:0.1}),
we find
\begin{equation}
U_{max} \approx 111, \qquad  \sup_{\Omega}\{|\omega(x,y)|)\}\approx 1.5 \times 10^{5}.
\label{eq:0.2}
\end{equation}
Together with $\nu=0.01$, and $A=1$ in Eq.~(\ref{eq:0.1}), this leads
to a Reynolds number $RE \approx 11,100$.

In the next section, we imagine the data in Figure 1 to be {\em desired or hypothetical data} at some suitable time $T > 0$, and explore the feasibility of reconstructing $\psi,~u,~v,~\omega,$ at $t=0,$ from these hypothetical values at time $T$.

\section{Uncertainty in backward in time Navier-Stokes reconstruction}
In dissipative evolution equations, reconstruction backward in time from imprecise data at some positive time $T$, is 
an ill-posed continuation problem. With appropriate stabilizing constraints, the uncertainty in such backward continuations is a function of how quickly
the given evolution equation forgets the
past as time advances. This is explored in such
publications as \cite{karen, knopslog, knops, payne, payne3,car00,car77,hao},  and the references therein. For the Navier-Stokes system in Eq.~(\ref{eq:0}),
the best-known results are due to Knops and Payne \cite{knops, payne}, and are based on logarithmic convexity arguments.

In Eq.~(\ref{eq:0}), let $\left \{\tilde{u}(x,y,t),~\tilde{v}(x,y,t),~\tilde{\omega}(x,y,t) \right \},$ be the solution on $\Omega \times [0,T)$ corresponding to approximate data
$\left \{\tilde{u}(\cdot,T),~\tilde{v}(\cdot,T),~\tilde{\omega}(\cdot,T) \right \},$ at a given positive time $~T$, and let $\left \{u(x,y,t),~v(x,y,t),~\omega(x,y,t) \right \}$ be the unique
solution associated with the true data $\left \{u(\cdot,T),~v(\cdot,T),~\omega(\cdot,T) \right \}$ at time $~T$.
Define $w(x,y,t)$ on $\Omega \times [0,T]$,  and $F(t)$ on $[0,T]$, as follows
\begin{equation}
\begin{array}{l}
(w(x,y,t))^2=\{u(x,y,t)-\tilde{u}(x,y,t)\}^2 + \{v(x,y,t)-\tilde{v}(x,y,t)\}^2, \\
\ \\
F(t)=\parallel w(\cdot, t) \parallel_2^2 = \int_{\Omega} (w(x,y,t))^2 dx dy.\\
\label{eq:1a}
\end{array}
\end{equation}
For a given small $\delta > 0$, assume that $|w(x,y,T)|$ is small enough that $F(T) \leq \delta~$.
What constraints must be placed on the solutions of Eq. (\ref{eq:0}) to ensure that $F(t)$ will remain small for $0 \leq t < T$ ? This is the
backward stability problem for the Navier-Stokes equations in the ${\cal{L}}^2(\Omega)$ norm.

Let $E,~Q,$ be prescribed positive constants. A velocity field $\{ u(x,y,t),~v(x,y,t) \}$ is said to belong to the class $\cal{E}$ if
\begin{equation}
\sup_{(x,t) \in \Omega \times [0,T]} \left\{(u^2+v^2) \right \} \leq E^2,
\label{eq:2}
\end{equation}
while it belongs to the class $\cal{Q}$ if
\begin{equation}
\sup_{(x,t) \in \Omega \times [0,T]} \left\{(u^2+v^2) + \omega^2 + (u_t^2+v_t^2) \right \} \leq Q^2.
\label{eq:3}
\end{equation}
Define $a,~b,~c,~\mu(t),~\Gamma(t),$ as follows,
\begin{equation}
\begin{array}{l}
a=2(E^2+1)/\nu,~~~~~b=Q^2(1+a/\nu),\\
\ \\
c=b/a,~~~~~~~~\mu(t)=(e^{at}-1)/(e^{aT}-1),~~~~0 \leq t \leq T, \\
\ \\
        \Gamma(t) \equiv \Gamma(E,Q,\nu,t,T)= \exp\{c(t-\mu(t)T)\},~~~~0 \leq t \leq T.
\label{eq:3a}
\end{array}
\end{equation}

If $\{ u(x,y,t),~v(x,y,t) \} \in \cal{E}$ and $\{\tilde{u}(x,y,t),~\tilde{v}(x,y,t) \} \in \cal{Q}$, it is shown in \cite{knops, payne} that then
\begin{equation}
F(t)F^{\prime \prime}(t)-(F^{\prime}(t))^2 \geq a F(t) F^{\prime}(t) - b F^2(t),
\label{eq:3b}
\end{equation}
from which it follows that with $\Gamma(t) \equiv \exp\{c(t-\mu(t)T)\}$,
\begin{equation}
        F(t) \leq \Gamma(t) \left\{ F(0) \right \}^{1-\mu(t)} \left\{ F(T) \right \}^{\mu(t)}, \qquad 0 \leq t \leq T.
\label{eq:4}
\end{equation}
A complete derivation for the result in Eq.~(\ref{eq:4}) may be found in \cite[Eqs.~9--12]{car6IPSE}.\\

Let  $M$ be an a-priori bound for $F(0)= \parallel w(\cdot, 0) \parallel_2^2$. With $F(T)\leq~\delta$, we then obtain
\begin{equation}
        \parallel w(\cdot, t) \parallel_2^2 = F(t) \leq \Gamma(t)  M^{1-\mu(t)} \delta^{\mu(t)}, \qquad 0 \leq t \leq T.
\label{eq:3c}
\end{equation}

However, as the Hurricane Ivan image in Figure 1 will show, even with very small $\delta > 0$, $\parallel w(\cdot, t) \parallel_2$ need not be small for $t < T$,
 because the factor $\Gamma(t) \equiv \exp\{c(t-\mu(t)T)\}$ in Eq.~(\ref{eq:3c})
may be extremely large, unless $T$ is extremely small.

Indeed, with $\nu=0.01$, and the values for $U_{max},~|\omega(x,y,0)|,$ given in Eq.~(\ref{eq:0.2}), we find
from Eqs.~(\ref{eq:2}, \ref{eq:3}, \ref{eq:3a}),
\begin{equation}
\begin{array}{l}
E^2 > U^2_{max} > 12,400,\qquad Q^2 > 2.19 \times 10^{10}, \\
\ \\
a \approx 2.48 \times 10^6,~~~~~b \approx  5.43 \times 10^{18},~~~~~c \approx 2.19 \times 10^{12}, \\
\ \\
        \mu(t)=(e^{at}-1)/(e^{aT}-1),\qquad \Gamma(t) = e^{\{c(t-\mu(t)T)\}}.
\label{eq:5}
\end{array} 
\end{equation}
With $T=1.0 \times 10^{-9},$ and $t=T/2$, we find $\mu(T/2)= 0.4996899859,$ and $\Gamma(T/2) \approx 1.97$. However, with $T=1.0 \times 10^{-8},$ and $t=T/2$, we now
find $\mu(T/2)= 0.4969000367,$ and $\Gamma(T/2) \approx 
3.05 \times 10^{29}$. Thus, in the Hurricane Ivan example in Section 2.1, if $T \geq 1.0 \times 10^{-8}$, the reconstruction uncertainty at $t=T/2$ given in Eq.~(\ref{eq:3c}), 
would appear to be too large for any realistic value of the data error
$~\delta~$ at time $T$.  We shall revisit that example in the instructive computational experiments described in Section 7 below.

\section{Stabilized leapfrog scheme for the 2D Navier-Stokes Equations}

Let $\Omega$ be the unit square in $R^2$, and consider the following 2D Navier-Stokes system for $(x,y) \in \Omega,$
\begin{equation}
\omega_t=L^{\dagger}\omega \equiv \nu \Delta \omega -u\omega_x -v\omega_y,~~~~~\Delta \psi=-\omega, \qquad 0 < t \leq T,\\
\label{eq:1.001}
\end{equation}
together with the homogeneous boundary conditions
\begin{equation}
\psi(x,y,t)=u(x,y,t)=v(x,y,t)=\omega(x,y,t)=0,~(x,y,t)\in \partial \Omega \times [0,T], 
\label{eq:1.0011}
\end{equation}
and the initial values
\begin{equation}
u(x,y,0)=u_0(x,y),~~~v(x,y,0)=v_0(x,y),~~~\omega(x,y,0)=\omega_0(x,y). 
\label{eq:1.002}
\end{equation}

The well-posed forward initial value problem
in Eq.~(\ref{eq:1.001}) becomes ill-posed if the time direction is reversed, and one wishes to recover $u(x,y,0)),~v(x,y,0),$ and $\omega(x,y,0),$
from given approximate values for $u(x,y,T),~v(x,y,T),~\omega(x,y,T).$ We contemplate both forward and time-reversed computations by allowing for possible {\em negative
time steps} $\Delta t$ in the leapfrog time-marching finite difference scheme described below.
With a given positive integer $N$, let $|\Delta t|=T/(N+1)$ be the time step magnitude, and let $\tilde{u}^n(x,y) \equiv \tilde{u}(x,y, n \Delta t), ~n=1,\cdots, N+1$,
denote the intended approximation to $u(x,y, n \Delta t)$, and likewise for $\tilde{v}^n(x,y)$, and $\tilde{\omega}^n(x,y).$ It is helpful to consider Fourier series expansions for
$\tilde{\omega}^n(x,y)$, on the unit square $\Omega$,
\begin{equation}
\tilde{\omega}^n(x,y)=\sum_{j,k = -\infty}^{\infty} \tilde{\omega}^n_{j,k} \exp\{2\pi i(jx+ky)\}, 
\label{eq:1.002A}
\end{equation}
with Fourier coefficients $\{\tilde{\omega}^n_{j,k}\}$ given by
\begin{equation}
\tilde{\omega}^n_{j,k}=\int_{\Omega} \tilde{\omega}^n(x,y) \exp\{-2\pi i (jx+ky)\} dxdy,
\label{eq:1.0022}
\end{equation}
With given fixed $\gamma > 0$ and $~p > 1$, define $\lambda_{j,k},~ \sigma_{j,k},$  as follows
\begin{equation}
\lambda_{j,k}=4 \pi^2 \nu (j^2+k^2), \qquad \sigma_{j,k}= \exp \{-\gamma |\Delta t| \lambda_{j,k}^p \}.
\label{eq:1.0023}
\end{equation}
For any $f(x,y) \in {\cal{L}}^2(\Omega)$, let $\{f_{j,k}\}$ be its Fourier coefficients as in Eq~(\ref{eq:1.0022}). Using
Eq.~(\ref{eq:1.0023}), define the linear operators $P$ and $S$ as follows
\begin{equation}
\begin{array}{l}
Pf=\sum_{j,k = -\infty}^{\infty} \lambda^p_{j,k} f_{j,k} \exp\{2\pi i(jx+ky)\}, \qquad \forall f \in {\cal{L}}^2(\Omega), \\
\ \\
Sf=\sum_{j,k = -\infty}^{\infty} \sigma_{j,k} f_{j,k} \exp\{2\pi i(jx+ky)\}, \qquad \forall f \in {\cal{L}}^2(\Omega).
\label{eq:1.004}
\end{array}
\end{equation}
As in \cite{carGEM, car1IPSE, car2IPSE, car3IPSE, Thermo, car5IPSE, car6IPSE, car7IPSE}, the operator $S$ is used as a stabilizing smoothing operator at each time step.
With the operator $L^{\dagger}$ as in Eq~(\ref{eq:1.001}), let $L^{\dagger}\tilde{\omega}^n \equiv \nu \Delta \tilde{\omega}^n -\tilde{u}^n \tilde{\omega}^n_x -\tilde{v}^n \tilde{\omega}^n_y$, and let
\begin{equation}
\tilde{\theta}^1=\omega_0(x,y), \qquad \tilde{\omega}^1=\omega_0(x,y) + \Delta t L^{\dagger}\omega_0.
\label{eq:1.000041}
\end{equation}
Consider the following stabilized leapfrog time-marching difference scheme for the system in Eq~(\ref{eq:1.001}), in which only the time variable is discretized, while the space variables remain continuous,
\begin{equation}
\begin{array}{l}
\Delta \tilde{\psi}^n=-\tilde{\omega}^n, \\
\ \\
\tilde{\theta}^{n+1}=S\tilde{\omega}^n, \\
\ \\
\tilde{\omega}^{n+1}=S \tilde{\theta}^n + 2 \Delta t S L^{\dagger}\tilde{\omega}^n,
 \qquad n=1,2,\cdots,N.
\label{eq:1.0041}
\end{array}
\end{equation}
The above semi-discrete problem is highly nonlinear. In the operator $L^{\dagger}$, the coefficients $\tilde{u}^n=\tilde{\psi}^n_y,~\tilde{v}^n=-\tilde{\psi}^n_x,$ are ultimately defined in terms of $~\tilde{\omega}^n$, which
is needed to obtain $~\tilde{\psi}^n$ by solving the Poisson problem $~\Delta \tilde{\psi}^n =- \tilde{\omega}^n$.
The analysis presented in Sections 5 and 6 below, although limited to a related simplified linear problem, is relevant to the above semi-discrete problem.
In Section 7, where
actual numerical computations are discussed, a uniform grid is imposed on $\Omega$, the space variables are discretized using centered differencing, and FFT algorithms are used to synthesize the smoothing operator $S$. A multigrid algorithm is used to solve the Poisson problem in Eq.~(\ref{eq:1.0041}).
\section{Fourier stability analysis in linearized problem}
Useful insight into the behavior of the nonlinear scheme in Eq.~(\ref{eq:1.0041}), can be gained
by analyzing a related linear problem with constant coefficients. We may naively view $L^{\dagger}$ in Eq.~(\ref{eq:1.002}), as a
linear operator with time dependent coefficients $u(x,y,t),~v(x,y,t),$ which need to be evaluated at every time step $t_n=n \Delta t$,
through an independent subsidiary calculation. If $|u(x,y,t)| \leq a,~|v((x,y,t)|\leq b,$ with positive constants $a,~b,$ we may study the numerical stability of the
system $~\tilde{\theta}^{n+1}=S \tilde{\omega}^n, ~~\tilde{\omega}^{n+1}=S\tilde{\theta}^n +2 \Delta t S (L^{\dagger} \tilde{\omega}^n)~$, in Eq.~(\ref{eq:1.0041}),
by considering the constant coefficient
linear operator $L\tilde{\omega}^n \equiv \nu \Delta \tilde{\omega}^n -a \tilde{\omega}^n_x -b \tilde{\omega}^n_y$, in lieu of $L^{\dagger} \tilde{\omega}^n$. Accordingly, we
examine the following evolution equation
\begin{equation}
\begin{array}{l}
\omega_t=L\omega \equiv \nu \Delta \omega -a\omega_x -b\omega_y,~~~0 < t \leq T, \\ 
\ \\
\omega(x,y,0)=\omega_0(x,y), 
\label{eq:3.001}
\end{array}
\end{equation}
together with homogeneous boundary conditions on $\partial \Omega$. With $\omega(x,y,t)$ the unique solution in Eq.~(\ref{eq:3.001}), let
\begin{equation} 
\theta^1=\omega(x,y,0),~~~~~~
\omega^n(x,y)=\omega(x,y, n\Delta t),~~ n \ge 1 .
\label{eq:3.0011}
\end{equation}
Then, the exact solution $\omega^n$ satisfies the following leapfrog system
\begin{equation}
\theta^{n+1}=\omega^n, ~~~\omega^{n+1}=\theta^n +2 \Delta t L \omega^n +\tau^n, \qquad n=1,2,\cdots,N,
\label{eq:3.0012}
\end{equation}
where the truncation error term $\tau^n$, is given by
\begin{equation}
\tau^n=\{(\Delta t)^3/6\} \{\omega_{ttt}(x,y,s)\}, ~~n|\Delta t| < |s| < (n+1)|\Delta t|,~~n=1,2,\cdots,N.
\label{eq:3.0013}
\end{equation}
Define the following two component vectors and matrix $G$
\begin{equation}
V^n=[\theta^n,~~\omega^n]^T,~~~~~~~ \Phi^n_\tau=[0,~~\tau_n]^T,
\label{eq:3.0014}
\end{equation}

\begin{equation}
G=
\begin{bmatrix}
0 & I \\
I & 2 \Delta t L \\
\end{bmatrix}
.
\label{eq:1.6661}
\end{equation}
One may then rewrite Eq.~(\ref{eq:3.0012}) as
\begin{equation}
V^{n+1}=GV^n+\Phi^n_{\tau}, \qquad n=1,2,\cdots,N.
\label{eq:3.0015}
\end{equation}
Define ${\cal{L}}^2(\Omega)$ norms
for two component vectors $W=[w^1~~w^2]^T$, and $2\times 2$ matrices $H$, as follows
\begin{equation}
\parallel W \parallel^2_2=(\parallel w^1 \parallel^2_2 + \parallel w^2 \parallel^2_2),~~~\parallel H \parallel_2=\sup_{\parallel W \parallel_2 \neq 0}\{\parallel HW \parallel_2/\parallel W \parallel_2 \}.
\label{eq:3.0016}
\end{equation}
Also, for functions $h(x,y,t)$ on $\Omega \times [0,T]$, define the norm $|||h|||_{2,\infty}$ as follows
\begin{equation}
|||h|||_{2,\infty}\equiv \sup_{0\leq t \leq T} \{\parallel h(\cdot,t)\parallel_2\}.
\label{eq:1.007}
\end{equation}
One may also write the stabilized leapfrog scheme for computing Eq.~(\ref{eq:3.001}) in matrix-vector notation.
Let
\begin{equation}
\tilde{\theta}^1(x,y)=\omega_0(x,y),~~~~
\tilde{\omega}^1(x,y)=\omega_0(x,y)+\Delta t L \omega_0.
\label{eq:3.0017}
\end{equation}
Then,
\begin{equation}
\begin{array}{l}
\omega(x,y,\Delta t)=\tilde{\omega}^1(x,y) +\tau^0, \\
\ \\
\tau^0=\{(\Delta t)^2/2\} \{\omega_{tt}(x,y,r)\}, ~~~0 < |r| < |\Delta t|.
\label{eq:3.0018}
\end{array}
\end{equation}
With $\tilde{V}^n=[\tilde{\theta}^n,~~\tilde{\omega}^n]^T$, and matrix $\Lambda$ defined as
\begin{equation}
\Lambda=
\begin{bmatrix}
S & 0 \\
0 & S \\
\end{bmatrix}
,
\label{eq:1.666}
\end{equation}
the stabilized marching leapfrog scheme for  Eq.~(\ref{eq:3.001}),
\begin{equation}
\tilde{\theta}^{n+1}= S\tilde{\omega}^n, ~~~~~\tilde{\omega}^{n+1}=S\tilde{\theta}^n +2\Delta t S L \tilde{\omega}^n, \qquad n=1,2, \cdots, N,
\label{eq:3.002}
\end{equation}
may be written as
\begin{equation}
\tilde{V}^{n+1}=\Lambda G \tilde{V}^n, \qquad n=1,2, \cdots, N.
\label{eq:3.0019}
\end{equation}
Note that with $\tau^0$ as in Eq.~(\ref{eq:3.0018}),
\begin{equation}
\parallel V^1-\tilde{V}^1\parallel_2=\parallel \tau^0 \parallel_2.
\label{eq:3.00191}
\end{equation}
Unlike the case in Eq.~(\ref{eq:1.0041}), the linear difference scheme in Eq.~(\ref{eq:3.0019}) is susceptible to Fourier analysis.
If $L \tilde{\omega}^n \equiv h(x,y)$, then
its Fourier coefficients $\{h_{j,k}\}$ are given by $h_{j,k}=g_{j,k}\tilde{\omega}^n_{j,k}$, where
\begin{equation}
g_{j,k}=-\{4 \pi^2 \nu (j^2+k^2) + 2 \pi i (aj+ bk)\}.
\label{eq:3.003}
\end{equation}
Hence, from Eqs.~(\ref{eq:1.0023}), (\ref{eq:1.004}) and (\ref{eq:3.002}),
\begin{equation}
\begin{array}{l}
\tilde{\theta}^{n+1}=\sum_{j,k =-\infty}^{\infty} ~\sigma_{j,k} \tilde{\omega}^{n}_{j,k} \exp \{2\pi i (jx+ky)\}, \\
\ \\
\tilde{\omega}^{n+1}= \sum_{j,k = -\infty}^{\infty}~ \sigma_{j,k} \{\tilde{\theta}^n_{j,k}+2\Delta t g_{j,k} \tilde{\omega}^n_{j,k} \} \exp\{2\pi i (jx+ky)\}.
\label{eq:3.004}
\end{array}
\end{equation}
\ \\
\ \\
Below, Lemmas 1 and 2, along with Theorems 1 and 2, are stated without proof. Using identical notation, proofs of these results may be found in \cite{car7IPSE}.

\newtheorem{lambdaQ}{Lemma}
\begin{lambdaQ}
Let $\lambda_{j,k},~\sigma_{j,k},$ be as in Eq.~(\ref{eq:1.0023}), and let $~g_{j,k}$ be as in Eq.~(\ref{eq:3.003}).
Choose a positive integer $J$ such that if $\lambda_J=4 \pi^2 \nu J$, we have
\begin{equation}
max_{(j^2+k^2) \leq J}~ \{|g_{j,k}|\} \leq 2 \lambda_J, \qquad |g_{j,k}| \leq 2 \lambda_{j,k},~\forall~(j^2+k^2) > J.
\label{eq:3.005}
\end{equation}
With $p > 1,$  choose $\gamma \geq 4(\lambda_J)^{1-p}$ in Eq.~(\ref{eq:1.0023}). Then,
\begin{equation}
\sigma_{j,k} \left(1+ 2 |\Delta t||g_{j,k}| \right) \leq 1 + 4|\Delta t| \lambda_J. \\
\label{eq:1.005}
\end{equation}
Hence,
\begin{equation}
\parallel \Lambda G \parallel_2 < 1+4|\Delta t| \lambda_J < \exp\{4 |\Delta t| \lambda_J\},
\label{eq:1.0051}
\end{equation}
and, for $n=1,2,\cdots,N$,
\begin{equation}
\parallel \tilde{V}^{n+1} \parallel_2  = \parallel (\Lambda G)^n \tilde{V}^1 \parallel_2 < \exp\{4n|\Delta t| \lambda_J\} \parallel \tilde{V}^1 \parallel_2.
\label{eq:1.0052}
\end{equation}
Therefore, with this choice of $(\gamma, p)$, the linear leapfrog scheme in Eq.~(\ref{eq:3.0019}), is unconditionally stable, marching forward or backward in time.
\end{lambdaQ}

\ \\
{\em Proof~}: See \cite[Lemma 5.1]{car7IPSE}
\newtheorem{Prelim}[lambdaQ]{Lemma}
\begin{Prelim}
Let $\omega^n(x,y) \equiv \omega(x,y, n \Delta t)$ be the exact solution in Eq.~(\ref{eq:3.001}). Let $\gamma,~ p,~ \lambda_{j,k},~ \sigma_{j,k},$ be as in Eq.~(\ref{eq:1.0023}).
Let $P$ and $S$ be as in Eq.~(\ref{eq:1.004}), let $L$ be the linear operator in Eq.~(\ref{eq:3.001}), and let $\tau^n$ be as in Eq.~(\ref{eq:3.0013}).
With the norm definition in Eq~(\ref{eq:1.007}), and $1 \leq n \leq N$,
\begin{eqnarray}
\parallel \tau^n \parallel_2 &\leq& (1/6)~|\Delta t|^3~|||\omega_{ttt} |||_{2,\infty}, \nonumber \\
\parallel \omega^n-S\omega^n \parallel_2 &\leq&  \gamma |\Delta t| ~||| P\omega |||_{2,\infty}, \nonumber \\
|\Delta t| \parallel L\omega^n-SL\omega^n \parallel_2 &\leq&  \gamma (\Delta t)^2~||| PL\omega |||_{2,\infty}.
\label{eq:12.006}
\end{eqnarray}
Moreover, with $V^n$ as in Eq.~(\ref{eq:3.0014}),
\begin{equation}
\parallel GV^n-\Lambda GV^n \parallel_2 \leq \gamma B \sqrt{3}~|\Delta t|~|||P\omega|||_{2,\infty},
\label{eq:12.0061}
\end{equation}
where the constant $B$ is given by
\begin{equation}
B=\{1+(8/3)(\Delta t)^2 |||PL\omega|||_{2,\infty}/|||P\omega|||_{2,\infty}\}^{1/2}.
\label{eq:12.0062}
\end{equation}
\end{Prelim}

\ \\
{\em Proof~}: See \cite[Lemma 5.2]{car7IPSE}
\ \\
\ \\
With the choice of $(\gamma, p)$ in Lemma 1,
the stabilized leapfrog scheme in Eq.~(\ref{eq:3.0019}) is unconditionally stable, but at the cost of introducing a small error at each
time step, whose cumulative effect will not vanish as $|\Delta t| \downarrow 0$. We now proceed to estimate the error $E^n=V^n - \tilde{V}^n,~~n=1,2,\cdots, N+1$.
Note that Theorems 1 and 2 below refer to the exact solution in Eq.~(\ref{eq:3.001}), which is a linear diffusion equation with constant coefficients.
Hence, $\omega(x,y,t)$ is
necessarily a very smooth function for $t > 0,$ in each of these Theorems. The same may not be true for the solution of the nonlinear problem in Eq.~(\ref{eq:1.001}). However,
in \cite{johnston,ghadimi}, highly accurate numerical schemes are advocated for the 2D Navier-Stokes equations.
Such schemes, which presume the existence of high order derivatives, are in common use.
Returning to Eq.~(\ref{eq:3.001}), the well-posed forward problem, with $\Delta t > 0$, is considered first.
\ \\
\newtheorem{Error}{Theorem}
\begin{Error}
With $\Delta t > 0$, let $\omega^{n}(x,y)$ be the unique exact solution of Eq.~(\ref{eq:3.001}) at $t=n \Delta t$.
With $\lambda_J$ be as in Lemma 1, let $p >1$ and let $\gamma=4\lambda_J^{1-p}$.
Let $V^n$ be as in Eq.~(\ref{eq:3.0015}),
let $\tilde{V}^n$ be the solution of the forward leapfrog scheme in Eq.~(\ref{eq:3.0019}),
and let $B$ be the constant defined in Eq.~(\ref{eq:12.0062}).
If $E^n=V^n-\tilde{V}^n$ denotes the error at $t_n= n \Delta t,~~n=1, 2, \cdots, N+1,$ we have, with the norm definitions in Eq.~(\ref{eq:1.007})
\begin{eqnarray}
\parallel E^n \parallel_2&\leq&\exp{(4\lambda_J t_n)}\parallel E^1 \parallel_2 +
	\sqrt{3}~B~(\lambda_J)^{-p}\{\exp{(4\lambda_J t_n)} -1\}~|||P\omega\|||_{2,\infty}\nonumber \\
~~~~&+&(24 \lambda_J)^{-1} (\Delta t)^2~\{\exp{(4\lambda_Jt_n)}-1\}~|||\omega_{ttt}|||_{2,\infty}. 
\label{eq:2.007}
\end{eqnarray}
Hence, as $\Delta t \downarrow 0$,
\begin{equation}
\parallel E^n \parallel_2 \leq  \exp{(4\lambda_J t_n)}\parallel E^1 \parallel_2 +\sqrt{3}~(\lambda_J)^{-p}~\{\exp{(4\lambda_J t_n)} -1\}~|||P\omega\|||_{2,\infty}+ O(\Delta t)^2.
\label{eq:2.00701}
\end{equation}
\end{Error}

\ \\
{\em Proof~}: See \cite[Theorem  5.1]{car7IPSE}
\ \\
\ \\
{\bf Remark 1}.
The first term on the right hand side in Eq.~(\ref{eq:2.00701}) allows for the possibility of erroneous initial data, leading to a positive
value for $\parallel E^1 \parallel_2$. The second term is the {\em stabilization penalty}, which is the price that must be paid for using the leapfrog scheme in the computation
of the well-posed forward linear problem in Eq.~(\ref{eq:3.001}). Without stabilization, the leapfrog scheme is unconditionally unstable for that forward problem, and no
error bound would be possible.
As is evident from the discussion in Section 3, small values of $T$ must necessarily
be contemplated in Navier-Stokes problems. For an important class of such problems with small $T > 0,$ and large $\lambda_J$, the value
of $(\lambda_J)^{-p}~\exp\{4\lambda_J T\}$, for some $p$ with $2.5  < p < 3.5$, may be sufficiently small to allow for useful computations. As an example,
with $T=1.0 \times 10^{-3},~ \lambda_J=2.0\times 10^{3},~p=3.25$, we find
$(\lambda_J)^{-p}~\exp\{4\lambda_J T\} < 5.572 \times 10^{-8}$.
\ \\

In the ill-posed backward problem, we contemplate marching backward in time from $t=T=(N+1)|\Delta t|$, using negative time time steps $\Delta t$. However, the needed initial data
$\omega_0(x,y)\equiv \omega(x,y,T)$, where $\omega(x,y,t)$ is the unique exact solution in Eq.~(\ref{eq:3.001}), may not be
known exactly. For this reason, the ill-posed backward problem is formulated as follows:

With $\Delta t < 0$, hypothetical  data $g(x,y)$ are given which are assumed to approximate the true data $~\omega(x, y, T)~$, while $~g + \Delta t Lg~$ approximates true data $~\omega(x,y, T-|\Delta t|)~$, both in the ${\cal{L}}^2(\Omega)$
norm. Moreover, the unique true backward solution $\omega(x,y,t)$ corresponding to the unknown exact data $\omega(x,y,T)$, is such that both $\omega(x,y,0)$ and $\omega(x,y, |\Delta t|)$
satisfy ${\cal{L}}^2(\Omega)$ bounds.
Specifically, for some unknown $\delta> 0$, and some $M  \gg~  \delta $,
\begin{equation}
\begin{array}{l}
\parallel g - \omega(\cdot, T) \parallel^2_2 + \parallel(g+\Delta t Lg)-\omega(\cdot,T-|\Delta t|)\parallel^2_2 \leq \delta^2, \\
\ \\
\parallel \omega(\cdot, 0) \parallel^2_2 + \parallel \omega (\cdot, |\Delta t|) \parallel^2_2 \leq M^2.
\label{eq:2.0083}
\end{array}
\end{equation}
Analogously to the well-posed forward problem in Eq.~(\ref{eq:3.0019}), we now choose $\Delta t < 0$ and consider the stabilized leapfrog scheme
marching backward from $t=T$ with the given hypothetical data $g(x,y)$ in Eq.~(\ref{eq:2.0083})
\begin{equation}
\tilde{V}^{n+1}=\Lambda G\tilde{V}^n,~~~~n=1,2,\cdots, N, ~~~~~~\tilde{V}^1=[g,~(g+\Delta t Lg)]^T,
\label{eq:2.0084}
\end{equation}
where, with $\delta,~M$ as in Eq.~(\ref{eq:2.0083}), and $ \Phi^n_{\tau}$ as in Eq.~(\ref{eq:3.0015}), the true solution satisfies
\begin{equation} 
V^{n+1}=GV^n + \Phi^n_{\tau},~~~n=1,2,\cdots N, ~~~~~~\parallel V^1-\tilde{V}^1 \parallel_2 \leq \delta,~~~~\parallel V^{N+1} \parallel_2 \leq M.
\label{eq:2.0085}
\end{equation}
\ \\
\newtheorem{ErrorB}[Error]{Theorem}
\begin{ErrorB}
With $g(x,y),~M,~\delta,$ as in Eq.~(\ref{eq:2.0083}), and $\Delta t < 0,$ let $V^{n}$ be the exact solution of the backward
problem in Eq.~(\ref{eq:2.0085}) at time $T-n|\Delta t|,~ n=1,2,\cdots, N+1.$
With  $p, ~\lambda_J,~$ as in Lemma 1, let $\gamma = 4(\lambda_J)^{1-p}$ in the smoothing operator $S$ in Eq.~(\ref{eq:1.0023}). Let $\tilde{V}^{n}$
be the corresponding solution of the stabilized backward leapfrog scheme in Eq.~(\ref{eq:2.0084}), let $B$ be the constant in Eq.~(\ref{eq:12.0062}),
and let $E^n\equiv V^{n}-\tilde{V}^{n}$
denote the error at time $T-n|\Delta t|$. Then, for $n=1,2,\cdots, N,$ with the norm definitions in Eq.~(\ref{eq:1.007})
\begin{eqnarray}
\parallel E^{n+1} \parallel_2 &\leq&
\sqrt{3}~B~(\lambda_J)^{-p}\left\{\exp(4 n|\Delta t|\lambda_J) -1 \right \} |||P\omega|||_{2,\infty} \nonumber \\
&+& (24 \lambda_J)^{-1}~(\Delta t)^2 \left\{\exp(4 n|\Delta t|\lambda_J) -1\right \} |||\omega_{ttt}|||_{2,\infty} \nonumber \\
&+& \delta \exp\{4 n|\Delta t|\lambda_J\}.
\label{eq:2.0087}
\end{eqnarray}
Hence, as $|\Delta t| \downarrow 0$,
\begin{eqnarray}
\parallel E^{n+1} \parallel_2 &\leq & 
\sqrt{3}(\lambda_J)^{-p} \left\{\exp(4n|\Delta t|\lambda_J) -1 \right \} |||P\omega|||_{2,\infty} \nonumber \\
&+& \delta \exp\{4 n|\Delta t|\lambda_J\}~+~O(\Delta t^2), \qquad n=1,2,\cdots, N. 
\label{eq:2.00871}
\end{eqnarray}
\end{ErrorB}

\ \\
{\em Proof~}: See \cite[Theorem  5.2]{car7IPSE}
\ \\
\ \\
{\bf Remark 2.} As was the case in Eq.~(\ref{eq:2.00701}), there is an additional error term, $\delta \exp\{4 n|\Delta t|\lambda_J\}$, 
in Eq.~(\ref{eq:2.00871}).
That term represents the necessary uncertainty in backward reconstruction from possibly erroneous data $g(x,y)$ at time $T$ in Eq.~(\ref{eq:3.001}).
\ \\
\subsection{Application to data assimilation}
In Theorems 1 and 2 above, define the constants $K_1$ through $K_4$ as follows, and consider the values shown in Table 1 below.\\
\begin{eqnarray}
 &~&~K_1=e^{4\lambda_J T},~K_2=\sqrt{3}~B(\lambda_J)^{-p} (K_1-1),~K_3=(24 \lambda_J)^{-1}(\Delta t)^2(K_1-1), \nonumber \\
 &~&~K_4=K_2 |||P\omega|||_{2,\infty} + K_3 |||\omega_{ttt}|||_{2,\infty}.
\label{eq:3.055}
\end{eqnarray}
\begin {center}
\ \\
{\em TABLE 1} \\
{\em Values of $K_1$ through $K_3$ in  Eq.~(\ref{eq:3.055}), with  following  parameter values:}\\
{$T=1\times10^{-5},~|\Delta t|=5\times10^{-8},~p=3.25,~\lambda_J=3.8\times10^{4},~(\lambda_J)^{-p}=1.3\times10^{-15}$.} \\
\ \\
\begin{tabular}{|c|c|c|} \hline
	{$K_1= e^{4\lambda_J T}$} & {$K_2=\sqrt{3}~B(\lambda_J)^{-p}(K_1-1)$} & {$K_3=(24 \lambda_J)^{-1}(\Delta t)^2(K_1-1)$} \\ \hline
        {$K_1 < 4.6$} & {$K_2 < 8.2 \times 10^{-15}$} & {$K_3 < 1.1 \times 10^{-20}$} \\ \hline
\end{tabular}
\end{center}
\ \\
As outlined in the Introduction, data assimilation applied to the system in Eq.~(\ref{eq:1.001}), is the problem of finding initial values $[u(.,0), v(.,0), \omega(.,0)]$, at $t=0$,
that can evolve into useful approximations to the given hypothetical data  at time $T$, as defined by  $\tilde{V}^1$ in Eq.~(\ref{eq:2.0084}) in the present linearized
problem.
If the true solution in  Eq.~(\ref{eq:3.001}) does not have exceedingly large values for 
$|||P\omega|||_{2,\infty}$ and $|||\omega_{ttt}|||_{2,\infty}$, 
the parameter values chosen in Table 1, together with Theorem 2, indicate that marching backward to time $t=0$ from the hypothetical data 
$\tilde{V}^1$ at time $T$, leads to an error $E(0)$ at $t=0$, satisfying
\begin{equation}
        \parallel E(0) \parallel_2 \leq \delta K_1 + K_4,
        \label{eq:3.056}
\end{equation}
with the constant $K_4$ defined in Eq.~(\ref{eq:3.055}) presumed small. Next, from Theorem 1, marching forward to time $T$, using the inexact computed 
initial values, leads to an error $E(T)$ at time $T$, satisfying
\begin{equation}
	\parallel E(T) \parallel_2 \leq K_1(\delta K_1 + K_4) + K_4.
        \label{eq:3.057}
\end{equation}
The error $E(T)$ in Theorem 1 is the difference at time $T$, between the unknown unique solution in Eq.~(\ref{eq:3.001}), and the computed numerical approximation to 
it provided by the stabilized forward explicit scheme. However, that unknown solution at time $T$ differs from the given hypothetical data at time $T$,
by an unknown amount $\delta $ in the $L^2$ norm. Hence the total error is
\begin{equation}
        \parallel E(T) \parallel_2 + \delta \leq \delta (1+K_1^2) +K_4 (1+K_1).  
        \label{eq:3.058}
\end{equation}
Therefore, data assimilation is successful only if the inexact computed initial values lead to a sufficiently small right hand side
in  Eq.(\ref{eq:3.058}). Clearly, the value of $\lambda_J T$, together with the unknown value of $\delta$, will  play a vital role.
As was the case in the nonlinear Knops-Payne uncertainty estimate following Eq.~(\ref{eq:5}) in Section 3,
increasing $T$ by a factor of ten may not be feasible,
even in the 
simplified linearized problem considered in Table 1. The value of the term $\delta (1+K_1^2)$ in Eq.(\ref{eq:3.058}), would increase from
$22~\delta$ to about $(1.6 \times 10^{13})~\delta$.

\section{Leapfrog nonlinear computational instability and the RAW filter}
The results in Theorems 1 and 2 indicate that the stabilizing approach in Eqs.~(\ref{eq:3.0019}) and (\ref{eq:2.0084}), is sound for the leapfrog scheme applied to the linearized
problem. 
However, for the $O(\Delta t)^2$ leapfrog scheme, linear stability is necessary
but not always sufficient in the presence of nonlinearities.
Leapfrog centered time differencing has been a mainstay of geophysical fluid dynamics computations. In well-posed linear wave propagation problems, a Courant condition is required for leapfrog computational stability. However, even when that Courant
condition is obeyed, serious instability can develop in nonlinear problems. This phenomenon, first reported in \cite{phil}, was  subsequently explored and analyzed in
\cite{fornberg, sanz},
and the references therein. Following \cite{phil,fornberg,sanz}, effective post processing time domain filtering techniques, that can prevent nonlinear instability, were developed in \cite{ass, williams1, williams2, amez,hurl}.
Such filtering is applied at every time step, and consists of replacing the computed solution with a specific linear combination of computations at previous time steps. 
These filtering techniques can also be usefully applied in the present 2D Navier-Stokes leapfrog computational  problem.

We first describe the Robert-Asselin-Williams filter (RAW), \cite{williams1,williams2}, as it applies to the stabilized leapfrog scheme for the forward nonlinear problem with
$\Delta t > 0,$ previously considered in Eq.~(\ref{eq:1.0041}).
With $L\omega^n \equiv \nu \Delta \omega^n -u^n \omega^n_x -v^n \omega^n_y$, and
\begin{equation}
\theta^1=\omega_0(x,y), \qquad \omega^1=\omega_0(x,y) + \Delta t L \omega_0,
\label{eq:1.000041222}
\end{equation}
this forward problem is defined by
\begin{equation}
\begin{array}{l}
\Delta \psi^n=-\omega^n, \\
\ \\
\theta^{n+1}=S\omega^n, \\
\ \\
\omega^{n+1}=S \theta^n + 2 \Delta t S L \omega^n, \qquad n=1,2,\cdots,N.
\label{eq:1.0041222}
\end{array}
\end{equation}
For $n=1,2,\cdots,N$, the forward problem is now modified as follows.
With positive constants $\xi,~\eta$, where $\xi = 0.53,$ and $0.01 \leq \eta \leq 0.2$, and with $\overline{\omega^1}=\omega^1,~\overline{\theta^1}=\theta^1,$
\begin{equation}
\begin{array}{l}
\Delta {\psi^n}=-\overline{\omega^n}, \\
\ \\
\theta^{n+1}=S\overline{\omega^n}, \\
\ \\
\omega^{n+1}=S \overline{\theta^n} + 2 \Delta t S L \overline{\omega^n}, \ \\
\ \\
\overline{\theta^{n+1}}=\theta^{n+1}+0.5\xi\eta\left(\omega^{n+1}-2 \theta^{n+1} +\overline{\theta^n}\right),\ \\
\ \\
\overline{\omega^{n+1}}=\omega^{n+1}-0.5\eta(1-\xi)\left(\omega^{n+1}-2 \overline{\theta^{n+1}} + \overline{\theta^n}\right ), 
\label{eq:1.00412223}
\end{array}
\end{equation}
with the filtered arrays $\overline{\theta^{n+1}},~\overline{\omega^{n+1}},$ overwriting the unfiltered arrays $\theta^{n+1},~\omega^{n+1},$ at each time step.
With given data $g(x,y)$ for $\omega(x,y,T)$ satisfying Eq.~(\ref{eq:2.0083}), and $\theta^1=g(x,y),~\omega^1=g+\Delta t Lg$, with $\Delta t < 0,$ the backward
problem is modified in exactly the same way.
The recommended value $\xi=0.53$ in \cite{williams1} is designed to maintain the $O(\Delta t)^2$ accuracy in the leapfrog scheme. \\
\ \\
{\bf Remark 3.}
Prior information about the solution is essential in ill-posed inverse problems. While accurate values for $M$ and $\delta$ in Eq.~(\ref{eq:2.0083})
are seldom known, interactive adjustment of the parameter pair $(\gamma, p)$ in Eqs.~(\ref{eq:1.0023}, \ref{eq:1.004}) in the smoothing operator $S$, often leads to useful reconstructions.
This process is similar to the manual tuning of an FM station, or the manual focusing of binoculars, and likewise requires user recognition of a correct solution.
Beginning with small values of $\gamma$ and $p$, chosen so as not to oversmooth the solution, a small number of successive trial reconstructions
are performed. Because of the underlying {\em explicit} marching difference scheme, this can be accomplished in a relatively short time. The values of $(\gamma, p)$ are increased slowly
if instability is detected, and are likewise decreased slowly to increase sharpness, provided no instability results. As is the case with binoculars, useful results are
generally obtained after relatively few trials.
There may be several possible good solutions. Typical values of $(\gamma, p)$ lie in the ranges $~10^{-14} \leq \gamma \leq 10^{-7},~~
2.5 \leq p \leq 3.5$.
\ \\

\section{Nonlinear computational experiments}
The stabilized leapfrog scheme, together with RAW filtering, can be successfully applied to the full 2D nonlinear Navier-Stokes equations. The computation proceeds
along the lines outlined in Section 4. At each time step, a multigrid algorithm is applied to the Poisson problem $\Delta \psi = -\omega$, with
$\omega$ computed at the previous time step. Centered space finite differencing is then applied to $\psi(x,y)$ to obtain $u(x,y), v(x,y)$, and 
to obtain $\omega$ at the next time step by solving the diffusion equation.

In each of Figures 2 through 5 below, the leftmost column presents a single preselected stream function image, as the same hypothetical data
at three distinct {\em increasing} values of $T > 0$. The middle column presents the corresponding initial value image at $t=0$, computed by the backward
marching stabilized leapfrog scheme. Finally, the rightmost column, presents the forward evolution to time $T > 0$, of the  previously computed
initial value image immediately to the left.

Evidently, given the same hypothetical stream function data, assimilation becomes less successful with increasing preselected $T > 0$.
The first row in each of Figures 2 through 5, is considered a successful example of 2D Navier-Stokes data assimilation. For each first row, the
corresponding behavior of the ${\cal{L}}^2$ norms for each of $\psi(x,y), u(x,y), v(x,y)$ and $\omega(x,y)$, in the assimilation process, 
is described in Tables 2 through 5. An important phenomenon is immediately noticeable in each of these tables. The ${\cal{L}}^2$ norm of the
evolved stream function $\psi(x,y)$ at time $T$, in the rightmost image in each first row, is a close approximation to that in the hypothetical
stream function in the leftmost image in that same first row. However, the ${\cal{L}}^2$ norms of the derivatives of $\psi(x,y)$ in the
evolved image at time $T$ in each first row, are substantially smaller than the corresponding ${\cal{L}}^2$ norms in the leftmost hypothetical
data image. This can be explained by the fact that the evolved image at time $T$ is a {\em computed solution} to the 2D Navier-Stokes equations. 
Such a solution would
necessarily be a smooth approximation to hypothetical data that might be far from an actual solution at the same time $T > 0.$
Indeed, Figure 6 shows how the stream function data $\psi(x,y)$, in the first row hypothetical and evolved Hurricane Ivan images, 
actually compare.
The approximating evolved data on the right, are noticeably smoother than  the hypothetical data on the left.

Finally, it is noteworthy that successful assimilation from a $T$ value on the order of $10^{-4}$, was found possible in the Hurricane Ivan image, given the uncertainty estimates following Eq.~(\ref{eq:5}) in Section 3. These estimates indicated a necessary value
for $T$ on the order of $10^{-9}$. However, rigorous uncertainty estimates must necessarily contemplate worse case
error accumulation scenarios, and may be too pessimistic in individual cases. The other examples in Tables 2, 3 and 5, also involve significantly
larger $T$ values than might be anticipated.

\section{Concluding Remarks}
Along with \cite{Burg,Advec,Thermoassim}, the results in the present paper invite  useful scientific debate and comparisons,
as to whether equally good or better
results might be achieved using the computational methods described in \cite{cintra, howard, blum, qizhi, arcucci, Antil, Chong, Lund,auroux1,ou,auroux2,auroux3,pozo,gosse,gomez, xu, tomislava,camposvelho}. As an alternative computational approach, backward marching stabilized explicit schemes may also be helpful
in verifying the validity of suspected hallucinations in computations involving  machine learning.

\begin{figure}
        \centerline{\includegraphics[width=5.0in]{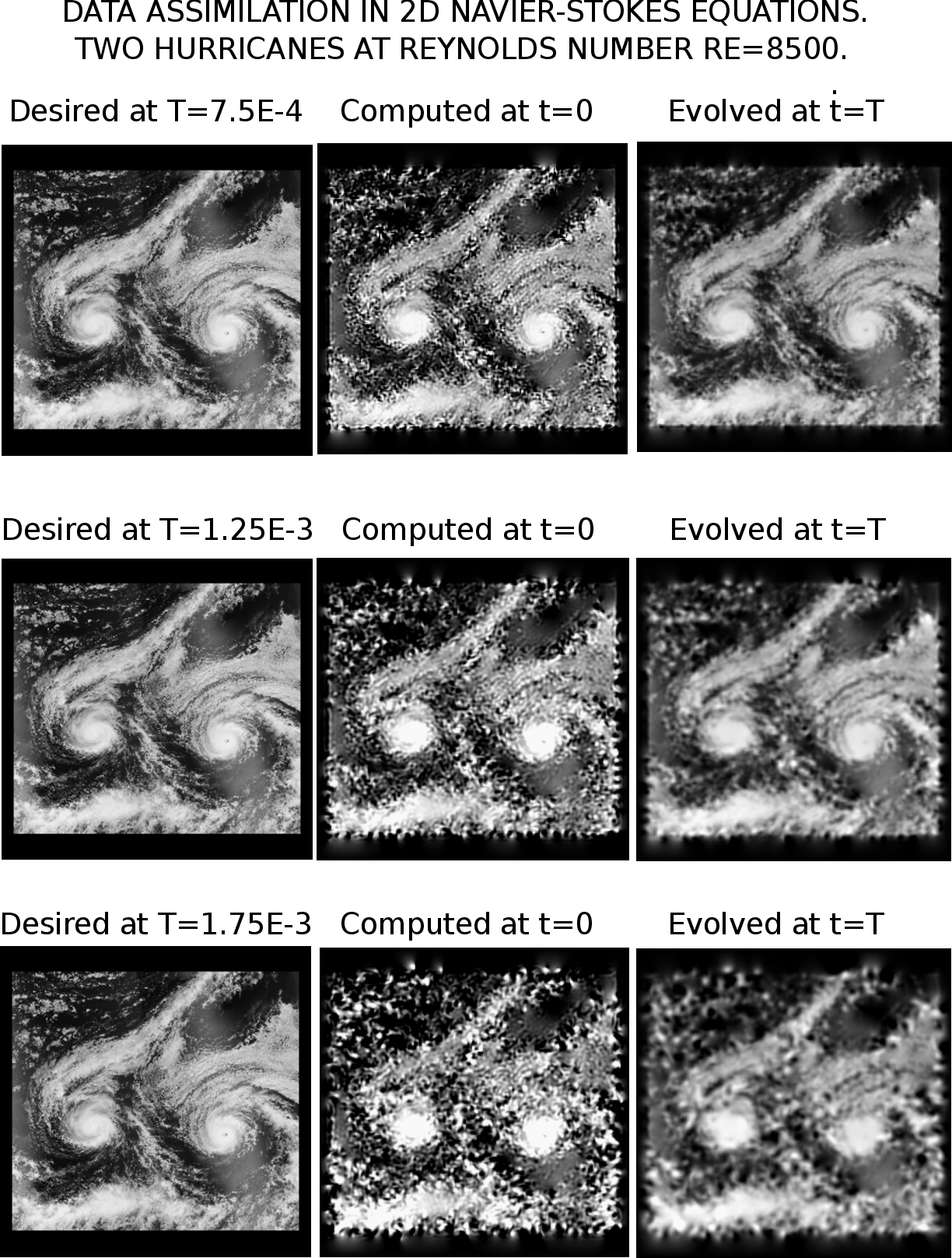}}
\caption{ As hypothetical data
	at various times $T > 0$, leftmost column features the August 2016 NASA Suomi NPP Satellite image of hurricanes Madeline and Lester approaching Hawaii. Middle column presents the corresponding initial values, obtained by solving the 2D Navier-Stokes equations backward in time.
	Rightmost column presents the forward evolution to time $T$ of these initial values, obtained by solving the same equations forward in time. Assimilation is
	successful if rightmost image is a useful approximation to leftmost image.}
\end{figure}

\begin{figure}
	\centerline{\includegraphics[width=5.0in]{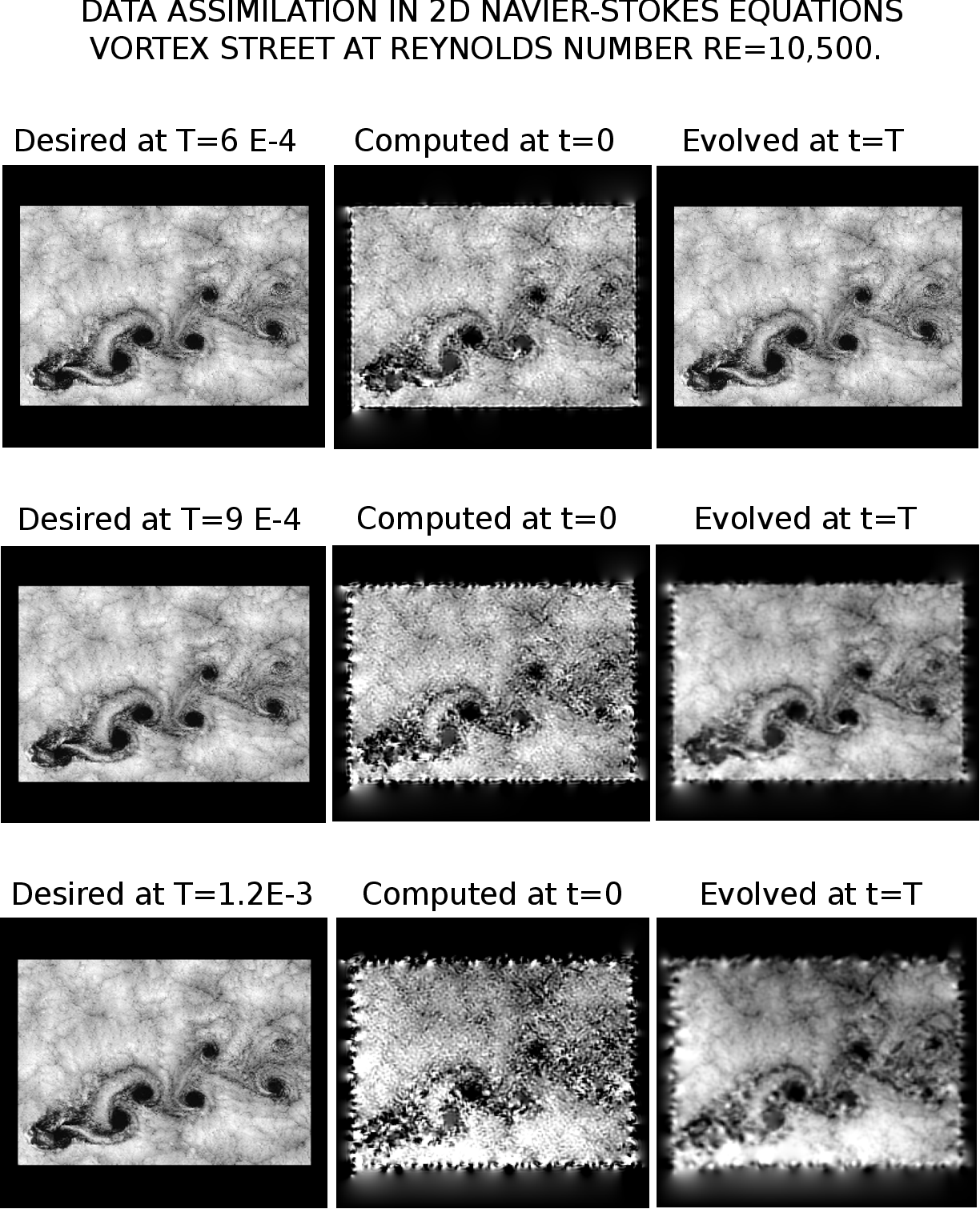}}
	\caption{As hypothetical data
        at various times $T > 0$, leftmost column features the September 1999 NASA Landsat 7 image of a von Karman vortex street in clouds off 
the Juan Fernandez islands. Middle column presents corresponding initial values, obtained by solving the 2D Navier-Stokes equations backward in time.
        Rightmost column presents the forward evolution to time $T$ of these initial values, obtained by solving the same equations forward in time. Assimilation is
        successful if rightmost image is a useful approximation to leftmost image.}
\end{figure}

\begin{figure}
	\centerline{\includegraphics[width=5.0in]{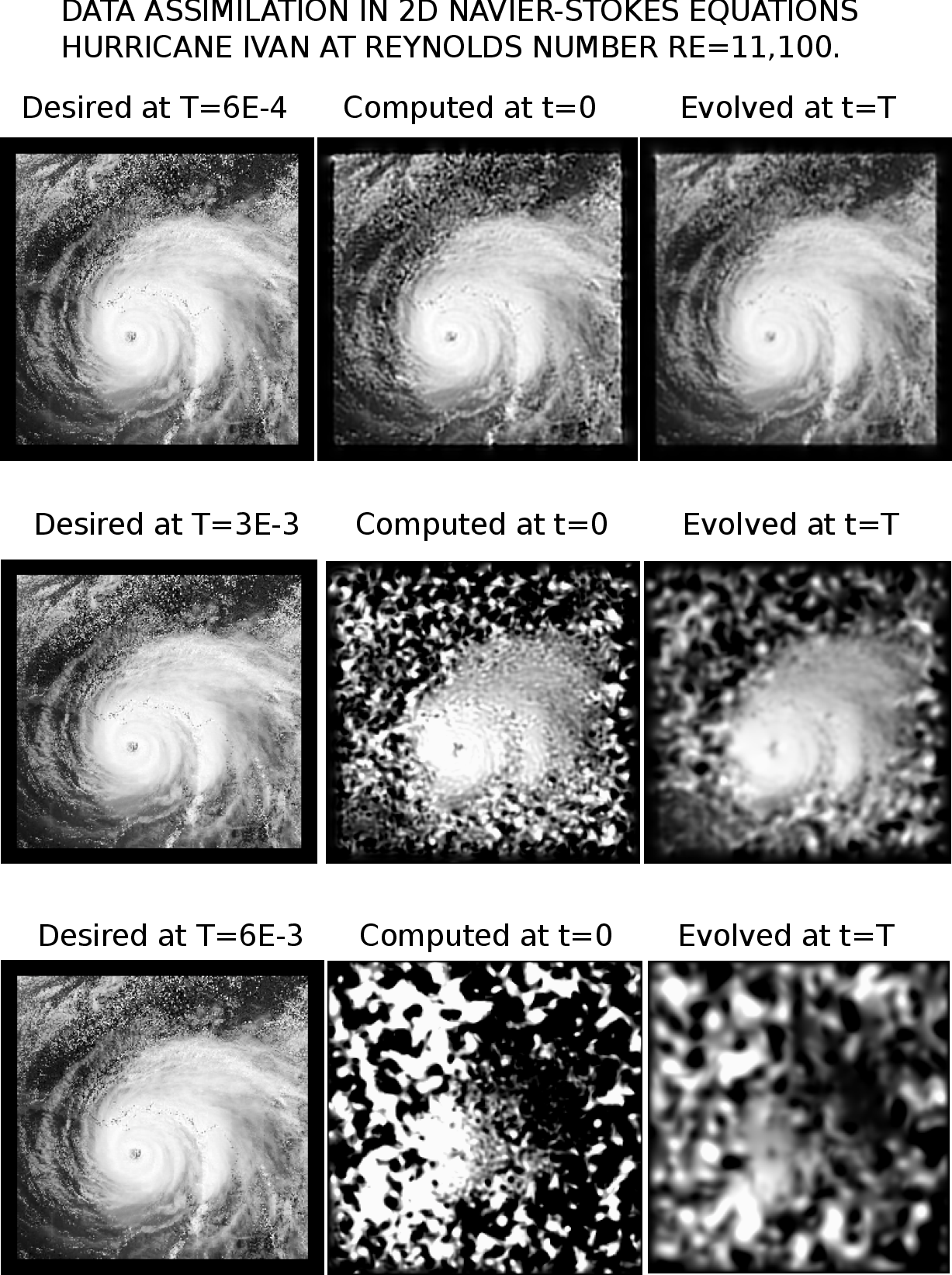}}
	\caption{As hypothetical data
	at various times $T > 0$, leftmost column features the September 2004 NASA Aqua Satellite image of Hurricane Ivan off Florida coast. 
	Middle column presents the corresponding initial values, obtained by solving the 2D Navier-Stokes equations backward in time.
        Rightmost column presents the forward evolution to time $T$ of these initial values, obtained by solving the same equations forward in time. Assimilation is
        successful if rightmost image is a useful approximation to leftmost image.}
\end{figure}
\begin{figure}
	\centerline{\includegraphics[width=5.0in]{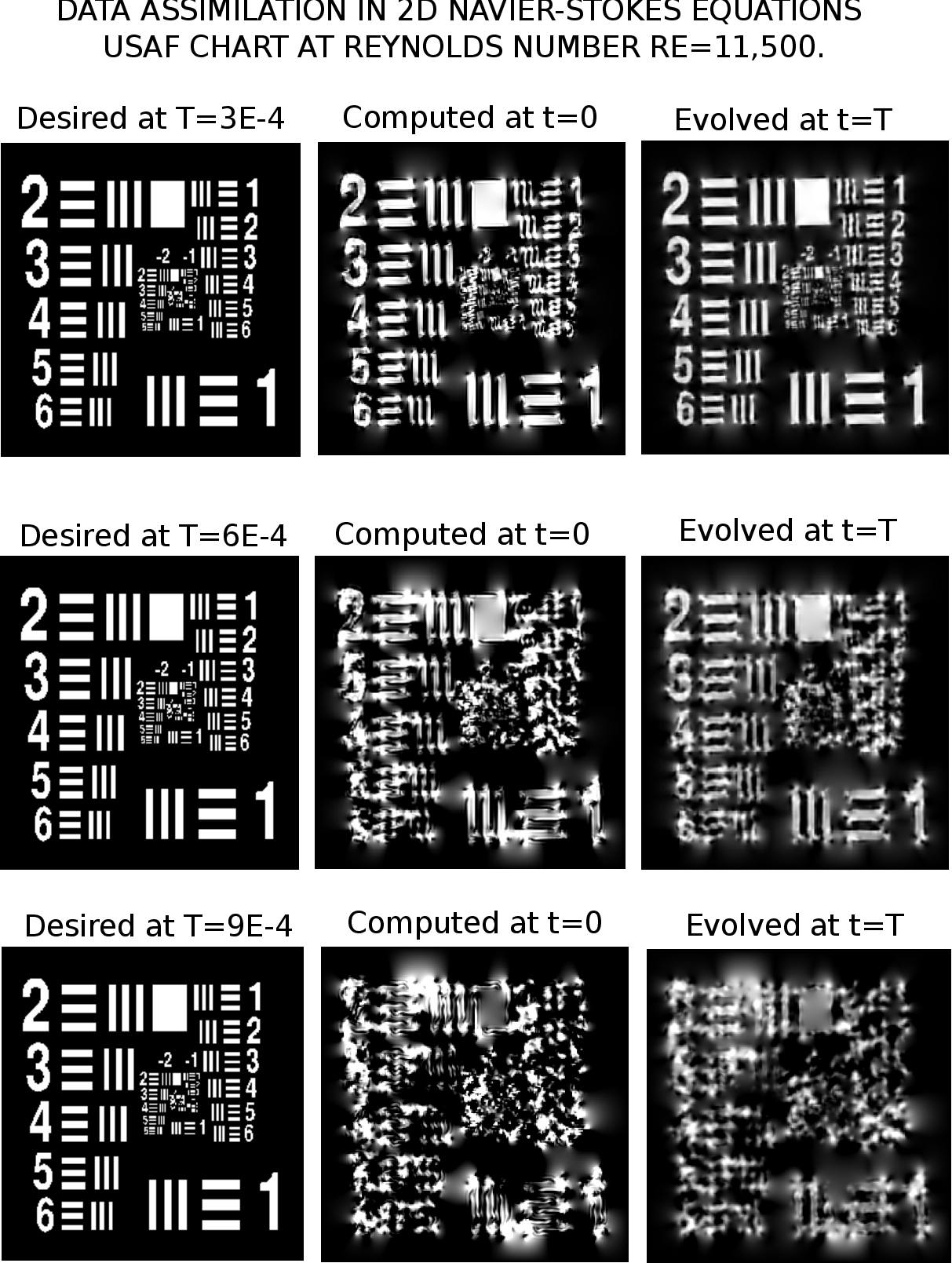}}
	\caption{As hypothetical data at various times $T > 0,$ leftmost column features the 1951 USAF Resolution Chart.
	Middle column presents the corresponding initial values, obtained by solving the 2D Navier-Stokes equations backward in time.
        Rightmost column presents the forward evolution to time $T$ of these initial values, obtained by solving the same equations forward in time. Assimilation is
        successful if rightmost image is a useful approximation to leftmost image.}
\end{figure}
\begin{figure}
\centerline{\includegraphics[width=4.75in]{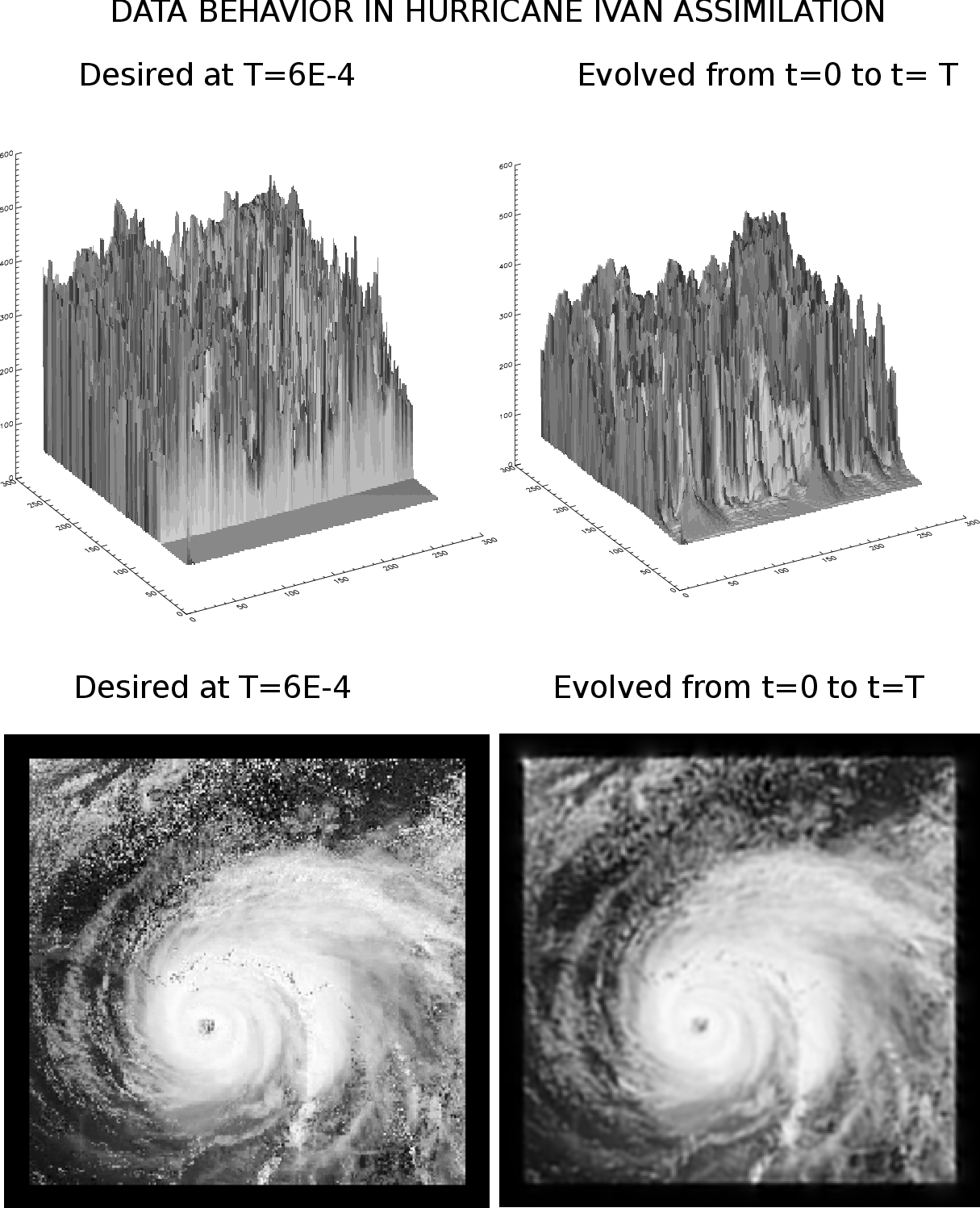}}
	\caption{The evolved Hurricane Ivan data and image in the right column, result from a computed solution of the Navier-Stokes 
	equations forward in time,
	and are a smooth approximation to the hypothetical data and image in the left column}
\end{figure}
\ \\
\ \\
\newpage
\begin {center}
{\em TABLE 2} \\
{\em Two Hurricanes image at $RE=8500$} \\
{\em ${\cal{L}}^2$-norm behavior in data assimilation from hypothetical data at $T=7.5E-4$} \\
\ \\
\begin{tabular}{|c|c|c|c|c|} \hline
	{\em Variable} & {\em Desired at $T$ norm}  & {\em Computed at $0$ norm}& {\em Evolved at $T$ norm} \\ \hline
	$\psi$ & 0.294 & 0.308 & 0.288    \\ \hline
$u=\psi_y$ & 12.3 & 20.6 & 8.5  \\ \hline
$v=-\psi_x$ & 12.7 & 20.4  & 8.8 \\ \hline
$\omega=-\Delta \psi$ & 14301 & 13986 & 3951 \\ \hline
\end{tabular}
\end{center}
\ \\
\begin {center}
{\em TABLE 3} \\
{\em Vortex Street image at RE=10,500}\\
{\em ${\cal{L}}^2$-norm behavior in data assimilation from hypothetical data at $T=6E-4$} \\
\ \\
\begin{tabular}{|c|c|c|c|c|} \hline
{\em Variable} & {\em Desired at $T$ norm} &{\em Computed at $0$ norm}  & {\em Evolved at $T$} norm \\ \hline
$\psi$ & 0.352 & 0.351 & 0.341 \\ \hline
$u=\psi_y$ & 10.0 & 14.1 & 7.3  \\ \hline
$v=-\psi_x$ & 10.7 & 15.0 & 7.9 \\ \hline
$\omega=-\Delta \psi$ & 10815 & 9579 & 3492 \\ \hline
\end{tabular}
\end{center}
\ \\
\begin {center}
{\em TABLE 4} \\
{\em Hurricane Ivan  image at RE=11,100} \\
{\em ${\cal{L}}^2$-norm behavior in data assimilation from hypothetical data at  $T=6E-4$} \\
\ \\
\begin{tabular}{|c|c|c|c|c|} \hline
{\em Variable} & {\em Desired at $T$ norm} &{\em Computed at $0$ norm }  & {\em Evolved at $T$ norm} \\ \hline
	$\psi$ & 0.356 & 0.362 & 0.350   \\ \hline
	$u=\psi_y$ & 13.32 & 18.21 & 8.46\\ \hline
$v=-\psi_x$ & 12.89 & 17.74 & 7.87 \\ \hline
$\omega=-\Delta \psi$ & 18094 &13207 & 4426 \\ \hline
\end{tabular}
\end{center}
\ \\
\begin {center}
{\em TABLE 5} \\
{\em USAF Chart image at RE=11,500} \\
{\em ${\cal{L}}^2$-norm behavior in data assimilation from hypothetical data at  $T=3E-4$} \\
\ \\
\begin{tabular}{|c|c|c|c|c|} \hline
{\em Variable} & {\em Desired at $T$ norm} &{\em Computed at $0$ norm}  & {\em Evolved at $T$ norm}  \\ \hline
$\psi$ & 0.263 & 0.300& 0.245  \\ \hline
	$u=\psi_y$ & 23.62 &  27.75 & 16.93 \\ \hline
$v=-\psi_x$ & 21.17 & 25.25 & 15.19 \\ \hline
	$\omega=-\Delta \psi$ & 15916 & 15110 & 7695 \\ \hline
\end{tabular}
\end{center}
\ \\


\newpage
\begin{thebibliography}{99}
	\bibitem{cintra} Cintra R, de Campos Velho HF, Cocke S.  Tracking the model: data assimilation by artificial neural network. 2016 International Joint
		Conference on Neural Networks (IJCNN), Vancouver, BC, Canada. 2016:403-410. DOI:10.1109/IJCNN.2016.7727227.
	\bibitem{howard}Howard LJ, Subramanian A, Hoteit I. A machine learning augmented data assimilation method for high resolution observations. Journal of Advances in Modeling Earth Systems. 2024;16,e2023MS003774.
		https://doi.org/10.1029/2023MS003774.
	\bibitem{blum} Blum J, Le Dimet F-X, Navon I-M. Data Assimilation for Geophysical Fluids. Handbook of Numerical Analysis. 2009;14:385-441. https://doi.org/10.1016/S1570-8659(08)00209-3 
		\bibitem{qizhi} He Q, Barajas-Solano D, Tartakovsky G, et al. Physics-informed neural networks for multiphysics data assimilation with
application to subsurface transport. Advances in Water Resources 2020; https://doi.org/10.1016/j.advwatres.2020.103610.
        \bibitem{arcucci} Arcucci R, Zhu J, Hu S, et al. Deep data assimilation: integrating deep learning with data assimilation.
Appl. Sci. 2021;11:1114. https://doi.org/10.3390/app11031114.

\bibitem{Antil} Antil H, Lohner R, Price R. Data assimilation with deep neural nets informed by nudging. https://arxiv.org/abs/2111.11505.  November 2021.
	\bibitem{Chong} Chen C, Dou Y, Chen J, et al. A novel neural network training framework with data assimilation. Journal of Supercomputing. 2022;78:19020--19045. https://doi.org/10.1007/s11227-04629-7
\bibitem{Lund}Lundvall J, Kozlov V, Weinerfelt P. Iterative methods for data assimilation for Burgers' equation. J. Inverse Ill-Posed Probl. 2006;14:505--535.
\bibitem{auroux1}Auroux D, Blum J. A nudging-based data assimilation method for oceanographuc problems: the back and forth nudging (BFN) algorithm. Proc. Geophys. 2008;15:305--319.
\bibitem{ou} Ou K, Jameson A. Unsteady adjoint method for the optimal control of advection and Burgers' equation using high order spectral difference method. 49th AIAA Aerospace Science Meeting, 4-7 January 2011. Orlando, Florida.
\bibitem{auroux2}Auroux D, Nodet M. The back and forth nudging algorithm for data assimilation problems: theoretical results on transport equations. ESAIM:COCV 2012;18:318--342.
\bibitem{auroux3}Auroux D, Bansart P, Blum J. An evolution of the back and forth nudging for geophysical data assimilation: application to Burgers equation and comparison.
Inverse Probl. Sci. Eng. 2013;21:399-419
\bibitem{pozo}Allahverdi N, Pozo A, Zuazua E. Numerical aspects of large-time optimal control of Burgers' equation. ESAIM Mathematical Modeling and Numerical Analysis 2016;50:1371--1401.
\bibitem{gosse}Gosse L, Zuazua E. Filtered gradient algorithms for inverse design problems of one-dimensional Burgers' equation. Innovative Algorithms and Analysis 2017;197--227.
\bibitem{camposvelho}de Campos Velho HF, Barbosa VCF, Cocke S. Special issue on inverse problems in geosciences. Inverse Probl. Sci. Eng. 2013;21:355-356. DOI:10.1080/17415977.2012.712532
	\bibitem{gomez}Gomez-Hernandez JJ, Xu T. Contaminant source identification in acquifers: a critical view. Math Geosci 2022;54:437--458.
\bibitem{xu} Xu T, Zhang W, Gomez-Hernandez JJ, et al. Non-point contaminant source identification in an acquifer using the ensemble smoother
with multiple data assimilation. Journal of Hydrology 2022; 606:127405.
\bibitem{tomislava} Vukicevic T, Steyskal M, Hecht M. Properties of advection algorithms in the context of variational data assimilation.
Monthly Weather Review 2001;129:1221--1231.
\bibitem{Burg}Carasso AS. Data assimilation in 2D viscous Burgers equation using a stabilized explicit finite difference scheme run backward in time. Inverse Probl. Sci. Eng. 2021;29:3475--3489. DOI:10.1080/17415977.2021.200947
        \bibitem{Advec}Carasso AS.  Data assimilation in 2D nonlinear advection diffusion equations, using an explicit stabilized leapfrog scheme run backwrd in time. NIST Technical Note 2227, July 12 2022.  DOI:10.6028/NIST.TN2227
	\bibitem{Thermoassim} Carasso AS. Data assimilation in 2D hyperbolic/parabolic systems using a stabilized explicit finite difference scheme run backward in time. Applied Mathematics in Science and Engineering, 2024;32:1,228641.
		DOI:10.1080/27690911.2023.228641
		\bibitem{richtmyer} Richtmyer RD, Morton KW. Difference Methods for Initial Value Problems. 2nd ed. New York
(NY): Wiley; 1967

\bibitem{carGEM}Carasso AS. Compensating operators and stable backward in time marching in nonlinear parabolic equations. Int J Geomath 2014;5:1--16.
\bibitem{car1IPSE}Carasso AS. Stable explicit time-marching in well-posed or ill-posed nonlinear parabolic equations. Inverse Probl. Sci. Eng. 2016;24:1364--1384.
\bibitem{car2IPSE}Carasso AS. Stable explicit marching scheme in ill-posed time-reversed viscous wave equations. Inverse Probl. Sci. Eng. 2016;24:1454--1474.
\bibitem{car3IPSE}Carasso AS. Stabilized Richardson leapfrog scheme in explicit stepwise computation of forward or backward nonlinear parabolic equations. Inverse Probl. Sci. Eng. 2017;25:1--24.
	        \bibitem{Thermo} Carasso AS. Stabilized  backward in time explicit marching schemes in the numerical computation of ill-posed time-reversed hyperbolic/parabolic systems.
Inverse Probl. Sci. Eng. 2018;1:1–32. DOI:10.1080/17415977.2018.1446952
\bibitem{car5IPSE}Carasso AS. Stable explicit stepwise marching scheme in ill- posed time-reversed 2D Burgers' equation. Inverse Probl. Sci. Eng. 2018;27(12):1-17. DOI:10.1080/17415977.2018.1523905

\bibitem{car6IPSE}Carasso AS. Computing ill-posed time-reversed 2D Navier-Stokes equations, using a stabilized explicit finite difference scheme marching backward in time.
Inverse Probl. Sci. Eng. 2019; DOI:10.1080/17415977.2019.1698564
\bibitem{car7IPSE}Carasso AS. Stabilized leapfrog scheme run backward in time, and the explicit $O(\Delta t)^2$ stepwise computation of ill-posed time-reversed 2D
Navier-Stokes equations. Inverse Probl. Sci. Eng. 2021; DOI:10.1080/17415977.2021.1972997
\bibitem{johnston}Johnston H, Liu JG. Finite difference schemes for incompressible flow based on local pressure boundary conditions. J. Comput. Phys. 2002;180:120--154.
\bibitem{ghadimi}Ghadimi P, Dashtimanesh A. Solution of 2D Navier-Stokes equation by coupled finite difference-dual reciprocity boundary element method. Applied Mathematical
Modelling 2011;35:2110--2121.
\bibitem{karen} Ames KA, Straughan B. Non-Standard and Improperly Posed Problems. New York (NY): Academic Press; 1997.
\bibitem{payne3} Payne LE. Some remarks on ill-posed problems for viscous fluids. Int. J. Engng Sci. 1992;30:1341--1347.
	\bibitem{knopslog} Knops RJ. Logarithmic convexity and other techniques applied to problems in continuum
mechanics. In: Knops RJ, editor. Symposium on non-well-posed problems and logarithmic convexity. Vol. 316, Lecture notes in mathematics. New York (NY): Springer-Verlag; 1973.
\bibitem{knops} Knops RJ, Payne LE. On the stability of solutions of the Navier-Stokes equations backward in time. Arch. Rat. Mech. Anal. 1968;29:331--335.
\bibitem{payne} Payne LE. Uniqueness and continuous dependence criteria for the Navier-Stokes equations. Rocky Mountain J. Math. 1971;2:641--660.
\bibitem{car00} Carasso AS. Reconstructing the past from imprecise knowledge of the present: Effective non-uniqueness in solving parabolic equations backward in time. Math. Methods Appl. Sci.
2012;36:249-–261.
\bibitem{car77}Carasso A. Computing small solutions of Burgers' equation backwards in time. J. Math. Anal. App. 1977;59:169--209.
\bibitem{hao}H\`{a}o DN, Nguyen VD, Nguyen VT. Stability estimates for Burgers-type equations backward in time. J. Inverse Ill Posed Probl. 2015;23:41--49.

\bibitem{phil} Phillips NA. An example of non-linear computational instability. In: Bolin B. editor, The Atmosphere and the Sea in Motion. Rockefeller Institute; New York (NY) 1959.
p. 501-504.
\bibitem{fornberg} Fornberg B. On the instability of Leapfrog and Crank-Nicolson approximations of a nonlinear partial differential equation. Math Comput. 1973;27:45-57.
\bibitem{sanz} Sanz-Serna JM. Studies in numerical nonlinear instability I. Why do leapfrog schemes go unstable ? SIAM J. Sci. Stat. Comput. 1985;6:923-938.
\bibitem{ass} Asselin R. Frequency filter for time integrations. Mon Weather Rev. 1972;100:487-490.
\bibitem{williams1} Williams PD. A proposed modification to the Robert-Asselin time filter. Mon Weather Rev. 2009;137:2538-2546.
\bibitem{williams2} Williams PD. The RAW filter: An improvement to the Robert-Asselin filter in semi-implicit integrations. Mon Weather Rev. 2011;139:1996-2007.
\bibitem{amez} Amezcua J, Kalnay E, Williams PD. The effects of the RAW filter on the climatology and forecast skill of the SPEEDY model. Mon Weather Rev. 2011;139:608-619.
\bibitem{hurl} Hurl N, Layton W, Li Y, et al. Stability analysis of the Crank-Nicolson-Leapfrog method with the Robert-Asselin-Williams time filter. BIT Numer Math. 2014;54:1009-1021.




\end{thebibliography}
\end{document}